\documentclass[a4paper]{amsart}
\usepackage[latin1]{inputenc}
\usepackage{latexsym}
\input amssym.def
\input amssym.tex
\newtheorem{definition}{Definition}[section]

\newtheorem{theorem}[definition]{Theorem}
\newtheorem{proposition}[definition]{Proposition}
\newtheorem{corollary}[definition]{Corollary}
\newtheorem{remark}[definition]{Remark}

\def\emph#1{{\bfseries\itshape{#1}}}
 1   
 1  

\def\R{\mathbb{R}}               

\def\lcf{\lbrack\! \lbrack}
\def\rcf{\rbrack\! \rbrack}

\newcommand\map[3]{#1\ \colon\ #2\longrightarrow#3}
\newcount\ancho
\newcount\anchom
\newcount\anchoa
\newcount\anchob
\newcommand{\ltilde}[2]{\ancho=#1 \anchom=\ancho \divide\anchom by 2
            \anchoa=\ancho \divide\anchoa by 4
        \anchob=\anchom \advance\anchob by \anchoa
      $\kern-5pt \begin{array}[b]{c}
                 \begin{picture}(1,1)(\anchom,0)
         \qbezier(0,2)(\anchoa,5)(\anchom,2)
         \qbezier(\anchom,2)(\anchob,-1)(\ancho,4)
         \qbezier(0,2)(\anchoa,4.5)(\anchom,1.8)
         \qbezier(\anchom,1.8)(\anchob,-1.5)(\ancho,4)
      \end{picture} \\[-4pt]
       \mbox{#2}
       \end{array} \kern-9pt$
       }
\begin{document}

\title{The Hamilton-Jacobi equation on Lie affgebroids}
\author[J.\ C.\ Marrero and D. Sosa]{J.C. Marrero and D. Sosa}
\address{J.C.\ Marrero and D. Sosa:
Departamento de Matem\'atica Fundamental, Facultad de
Ma\-te\-m\'a\-ti\-cas, Universidad de la Laguna, La Laguna,
Tenerife, Canary Islands, Spain} \email{jcmarrer@ull.es,
dnsosa@ull.es}

\thanks{This work has been partially supported by MICYT (Spain)
Grant BMF 2003-01319}



\begin{abstract}
\noindent The Hamilton-Jacobi equation for a Hamiltonian section
on a Lie affgebroid is introduced and some examples are discussed.
\end{abstract}

\maketitle

\setcounter{section}{0}

\section{Introduction}

Recently, in \cite{LMM} (see also \cite{M2}) the authors developed
a Hamiltonian description of Mechanics on Lie algebroids. If
$\tau_{E}:E\to M$ is a Lie algebroid on $M$ then, in this
description, the role of the cotangent bundle of the configuration
manifold is played by the $E$-tangent bundle ${\mathcal T}^EE^*$
to $E^*$ (the prolongation of $E$ over $\tau_{E}^*:E^*\to M$ in
the terminology of \cite{LMM}). One may construct the canonical
symplectic 2-section associated with the Lie algebroid $E$ as a
closed non-degenerate section $\Omega_E$ of the vector bundle
$\wedge^2({\mathcal T}^EE^*)^*\to E^*$. Then, given a Hamiltonian
function $H:E^*\to\R$, the Hamiltonian section associated with $H$
is the section $\xi_H$ of ${\mathcal T}^EE^*\to E^*$ characterized
by the equation $i_{\xi_H}\Omega_E=d^{{\mathcal T}^EE^*}H$,
$d^{{\mathcal T}^EE^*}$ being the differential of the Lie
algebroid ${\mathcal T}^EE^*\to E^*$. The integral sections of
$\xi_H$ are the solutions of the Hamilton equations for $H$. In
fact, these solutions are just the integral curves of the
Hamiltonian vector field of $H$ with respect to the linear Poisson
structure on $E^*$ associated with the Lie algebroid $E$.

Using the canonical symplectic section $\Omega_E$, one may also
give a description of the Hamil\-to\-nian Mechanics on the Lie
algebroid $E$ in terms of Lagrangian submanifolds of symplectic
Lie algebroids (see \cite{LMM}). An alternative approach, using
the linear Poisson structure on $E^*$ and the canonical
isomorphism between $T^*E$ and $T^*E^*$ was discussed in
\cite{GGU3}.

In \cite{LMM}, the authors also introduced the Hamilton-Jacobi
equation for a Hamiltonian function $H:E^*\to \R$ and they proved
that knowing one solution of the Hamilton-Jacobi equation
simplifies the search of trajectories for the corresponding
Hamiltonian vector field.

On the other hand, in \cite{GGrU,MMeS} a possible generalization
of the notion of a Lie algebroid to affine bundles is introduced
in order to create a geometric model which provides a natural
framework for a time-dependent version of Lagrange equations on
Lie algebroids (see also \cite{GGU2,M,SMM}). The resultant objects
are called Lie affgebroid structures (in the terminology of
\cite{GGrU}). If $\tau_A:A\to M$ is an affine bundle modelled on
the vector bundle $\tau_V:V\to M$, $\tau_{A^+}:A^+=Aff(A,\R)\to M$
is the dual bundle to $A$ and $\widetilde{A}=(A^+)^*$ is the
bidual bundle, then a Lie affgebroid structure on $A$ is
equivalent to a Lie algebroid structure on $\widetilde{A}$ such
that the distinguished section $1_A$ of $\tau_{A^+}:A^+\to M$
(corresponding to the constant function $1$ on $A$) is a 1-cocycle
in the Lie algebroid cohomology complex of $\widetilde{A}$. Now,
if $h:V^*\to A^+$ is a Hamiltonian section (that is, $h$ is a
section of the canonical projection $\mu:A^+\to V^*$) then one may
construct a cosymplectic structure $(\Omega_h,\eta)$ on the
prolongation of $\widetilde{A}$ over the fibration
$\tau_V^*:V^*\to M$ and one may consider the Reeb section $R_h$ of
$(\Omega_h,\eta)$. The integral sections of $R_h$ are just the
solutions of the Hamilton equations for $h$ (see \cite{IMPS,M}).
Alternatively, one may prove that the solutions of the Hamilton
equations for $h$ are the integral curves of the Hamiltonian
vector field of $h$ on $V^*$ with respect to the canonical
aff-Poisson structure on the line affine bundle $\mu:A^+\to V^*$
(see \cite{IMPS}). Aff-Poisson structures on line affine bundles
were introduced in \cite{GGrU} (see also \cite{GGU2}) as the
affine analogs of Poisson structures. The existence of an
aff-Poisson structure on the affine bundle $\mu:A^+\to V^*$ is a
consequence of some general results proved in \cite{GGU2}.

The aim of this paper is to introduce the Hamilton-Jacobi equation
for a Hamiltonian section $h$ and, then, to prove that knowing one
solution of the Hamilton-Jacobi equation simplifies the search of
trajectories for the corresponding Hamiltonian section (see
Theorem \ref{thHJ}).

The paper is organized as follows. In Section 2, we recall some
definitions and results about Lie algebroids and Lie affgebroids
which will be used in the rest of the paper. The Hamiltonian
formalism on Lie affgebroids is developed in Section 3. The
Hamilton-Jacobi equation for a Hamiltonian section on a Lie
affgebroid $A$ is introduced in Section $4$ and, then, we prove
the main result of the paper (see Theorem \ref{thHJ}). Some
examples are discussed in the last section (Section 5). In
particular, if $A$ is a Lie algebroid then we recover the
Hamilton-Jacobi equation considered in \cite{LMM} (see Example
\ref{ex5.0}). When $A$ is the 1-jet bundle of local sections of a
trivial fibration $\tau:\R\times P\to \R$ then the Hamiltonian
section $h$ may be considered as a time-dependent Hamiltonian
function $H:\R\times T^*P\to \R$ and the resultant equation is
just the classical time-dependent Hamilton-Jacobi equation (see
Example \ref{ex5.1}). Finally, we obtain some results about the
Hamilton-Jacobi equation on the so-called Atiyah affgebroid
associated with a principal G-bundle $p:Q\to M$ and a fibration
$\nu:M\to \R$ (see Example \ref{ex5.2}).

\section{Lie algebroids and Lie affgebroids}
\subsection{Lie algebroids}\label{Lie-algebroids}
Let $E$ be a vector bundle of rank $n$ over the manifold $M$ of
dimension $m$ and $\tau_{E}:E\rightarrow M$ be the vector bundle
projection. Denote by $\Gamma(\tau_{E})$ the
$C^{\infty}(M)$-module of sections of $\tau_{E}:E\rightarrow M$. A
{\it Lie algebroid structure} $(\lcf\cdot,\cdot\rcf_{E},\rho_{E})$
on $E$ is a Lie bracket $\lcf\cdot,\cdot\rcf_{E}$ on the space
$\Gamma(\tau_{E})$ and a bundle map $\rho_{E}:E\rightarrow TM$,
called {\it the anchor map}, such that if we also denote by
$\rho_{E}:\Gamma(\tau_{E})\rightarrow\frak{X}(M)$ the homomorphism
of $C^{\infty}(M)$-modules induced by the anchor map then $\lcf
X,fY\rcf_{E}=f\lcf X,Y\rcf_{E}+\rho_{E}(X)(f)Y,$ for
$X,Y\in\Gamma(\tau_{E})$ and $f\in C^{\infty}(M)$. The triple
$(E,\lcf\cdot,\cdot\rcf_{E},\rho_{E})$ is called a {\it Lie
algebroid over $M$} (see \cite{Ma}). In such a case, the anchor
map $\rho_{E}:\Gamma(\tau_{E})\rightarrow\frak{X}(M)$ is a
homomorphism between the Lie algebras
$(\Gamma(\tau_{E}),\lcf\cdot,\cdot\rcf_{E})$ and
$(\frak{X}(M),[\cdot,\cdot])$.

If $(E,\lcf\cdot,\cdot\rcf_{E},\rho_{E})$ is a Lie algebroid, one
may define a cohomology operator which is called {\it the
differential of $E$},
$d^E:\Gamma(\wedge^k\tau_{E}^*)\longrightarrow\Gamma(\wedge^{k+1}\tau_{E}^*)$,
as follows
\begin{equation}\label{dE}
\begin{array}{lcl}
(d^E\mu)(X_0,\dots,X_k)&=&\displaystyle\sum_{i=0}^k(-1)^i\rho_{E}(X_i)(\mu(X_0,\dots,\widehat{X_i},\dots,X_k))\\
 &+&\displaystyle\sum_{i<j}(-1)^{i+j}\mu(\lcf
 X_i,X_j\rcf_{E},X_0,\dots,\widehat{X_i},\dots,\widehat{X_j},\dots,X_k),
\end{array}
\end{equation}
for $\mu\in\Gamma(\wedge^k\tau_{E}^*)$ and
$X_0,\dots,X_k\in\Gamma(\tau_{E})$. Moreover, if
$X\in\Gamma(\tau_{E})$, one may introduce, in a natural way, {\it
the Lie derivate with respect to $X$}, as the operator ${\mathcal
L}_X^E:\Gamma(\wedge^k\tau_{E}^*)\longrightarrow\Gamma(\wedge^{k}\tau_{E}^*)$
given by ${\mathcal L}_X^E=i_X\circ d^E+d^E\circ i_X.$

If $E$ is the standard Lie algebroid $TM$ then the differential
$d^E=d^{TM}$ is the usual exterior differential associated with
$M$, which we will denote by $d_0$.

If we take local coordinates $(x^i)$ on $M$ and a local basis
$\{e_\alpha\}$ of sections of $E$, then we have the corresponding
local coordinates $(x^i,y^\alpha)$ on $E$, where $y^\alpha(a)$ is
the $\alpha$-th coordinate of $a\in E$ in the given basis. Such
coordinates determine local functions $\rho_\alpha^i$,
$C_{\alpha\beta}^{\gamma}$ on $M$ which contain the local
information of the Lie algebroid structure, and accordingly they
are called the {\it structure functions of the Lie algebroid.}
They are given by
\[
\rho_{E}(e_\alpha)=\rho_\alpha^i\frac{\partial }{\partial
x^i}\;\;\;\mbox{ and }\;\;\; \lcf
e_\alpha,e_\beta\rcf_{E}=C_{\alpha\beta}^\gamma e_\gamma.
\]
%
%

If $f\in C^\infty(M)$, we have that
$$d^E f=\frac{\partial f}{\partial x^i}\rho_\alpha^i e^\alpha,$$
where $\{e^\alpha\}$ is the dual basis of $\{e_\alpha\}$. On the
other hand, if $\theta\in \Gamma(\tau_{E}^*)$ and
$\theta=\theta_\gamma e^\gamma$ it follows that
$$d^E \theta=(\frac{\partial \theta_\gamma}{\partial
x^i}\rho^i_\beta-\frac{1}{2}\theta_\alpha
C^\alpha_{\beta\gamma})e^{\beta}\wedge e^\gamma.$$

For a Lie algebroid $(E,\lcf\cdot,\cdot\rcf_{E},\rho_{E})$ over
$M$ we may consider the generalized distribution ${\mathcal F}^E$
on $M$ whose characteristic space at point $x\in M$ is given by
${\mathcal F}^E(x)=\rho_{E}(E_x)$, where $E_x$ is the fibre of $E$
over $x$. The distribution ${\mathcal F}^E$ is finitely generated
and involutive. Thus, ${\mathcal F}^E$ defines a generalized
foliation on $M$ in the sense of Sussmann \cite{Su}. ${\mathcal
F}^E$ is the {\it Lie algebroid foliation on $M$} associated with
$E$.

Note that if $f\in C^{\infty}(M)$ then $d^Ef=0$ if and only if $f$
is constant on the leaves of ${\mathcal F}^E$.

Now, suppose that $(E,\lcf\cdot,\cdot\rcf_{E},\rho_{E})$ and
$(E',\lcf\cdot,\cdot\rcf_{E'},\rho_{E'})$ are Lie algebroids over
$M$ and $M'$, respectively, and that $F:E\to E'$ is a vector
bundle morphism over the map $f:M\to M'.$ Then $(F,f)$ is said to
be a {\it Lie algebroid morphism} if
$$d^{E}((F,f)^*\phi')=(F,f)^*(d^{E'}\phi'),\mbox{ for }\phi'\in\Gamma(\wedge^k(\tau_{E'})^*) \mbox{ and for all } k.$$
Note that $(F,f)^*\phi'$ is the section of the vector bundle
$\wedge^k E^*\rightarrow M$ defined by
$$((F,f)^*\phi')_{x}(a_1,\dots,a_k)=\phi'_{f(x)}(F(a_1),\dots,F(a_k)),$$
for $x\in M$ and $a_1,\dots,a_k\in E_{x}$. If $(F,f)$ is a Lie
algebroid morphism, $f$ is an injective immersion and
$F_{|E_x}:E_x\rightarrow E'_{f(x)}$ is injective, for all $x\in
M$, then $(E,\lcf\cdot,\cdot\rcf_{E},\rho_{E})$ is said to be a
{\it Lie subalgebroid} of
$(E',\lcf\cdot,\cdot\rcf_{E'},\rho_{E'})$.

\subsubsection{The prolongation of a Lie algebroid over a
fibration}\label{sec1.1.1}

In this section, we will recall the definition of the Lie
algebroid structure on the prolongation of a Lie algebroid over a
fibration (see \cite{HM,LMM}).

Let $(E,\lcf\cdot,\cdot\rcf_{E},\rho_{E})$ be a Lie algebroid of
rank $n$ over a manifold $M$ of dimension $m$ with vector bundle
projection $\tau_{E}:E\rightarrow M$ and $\pi:M'\rightarrow M$ be
a fibration.

We consider the subset ${\mathcal T}^EM'$ of $E\times TM'$ and the
map $\tau_{E}^{\pi}:{\mathcal T}^EM'\rightarrow M'$ defined by
$${\mathcal T}^EM'=\{(b,v')\in E\times
TM'/\rho_{E}(b)=(T\pi)(v')\},\;\;\;\;\tau_{E}^{\pi}(b,v')=\pi_{M'}(v'),$$
where $T\pi:TM'\rightarrow TM$ is the tangent map to $\pi$ and
$\pi_{M'}:TM'\rightarrow M'$ is  the canonical projection. Then,
$\tau_{E}^{\pi}:{\mathcal T}^EM'\rightarrow M'$ is a vector bundle
over $M'$ of rank $n+\dim M'-m$ which admits a Lie algebroid
structure $(\lcf\cdot,\cdot\rcf_{E}^{\pi},\rho_{E}^{\pi})$
characterized by
$$\lcf(X\circ\pi,U'),(Y\circ\pi,V')\rcf_{E}^{\pi}=(\lcf X,Y\rcf_{E}\circ\pi
,[U',V']),\;\;\rho_{E}^{\pi}(X\circ\pi,U')=U',$$
for all  $X,Y\in \Gamma(\tau_{E})$  and $U', V'$ vector fields
which are $\pi$-projectable to $\rho_{E}(X)$ and $\rho_{E}(Y)$,
res\-pec\-tively.
 $({\mathcal
T}^EM',\lcf\cdot,\cdot\rcf_{E}^{\pi},\rho_{E}^{\pi})$ is called
{\it the prolongation of the Lie algebroid $E$ over the fibration
$\pi$ or the $E$-tangent bundle to $M'$} (for more details, see
\cite{HM,LMM}).

Next, we consider a particular case of the above construction. Let
$E$ be a Lie algebroid over a manifold $M$ with vector bundle
projection $\tau_{E}:E\rightarrow M$ and ${\mathcal T}^EE^*$ be
the prolongation of $E$ over the projection
$\tau_{E}^*:E^*\rightarrow M$. ${\mathcal T}^EE^*$ is a Lie
algebroid over $E^*$ and we can define a canonical section
$\lambda_E$ of the vector bundle $({\mathcal T}^EE^*)^*\rightarrow
E^*$ as follows. If $a^*\in E^*$ and $(b,v)\in({\mathcal
T}^EE^*)_{a^*}$ then
\begin{equation}\label{lambdaE}
\lambda_E(a^*)(b,v)=a^*(b).
\end{equation}
$\lambda_E$ is called the {\it Liouville section} associated with
the Lie algebroid $E.$

Now, one may consider the nondegenerate $2$-section
$\Omega_E=-d^{{\mathcal T}^EE^*}\lambda_E$ of ${\mathcal
T}^EE^*\rightarrow E^*$. It is clear that $d^{{\mathcal
T}^EE^*}\Omega_E=0$. In other words, $\Omega_E$ is a symplectic
section. $\Omega_E$ is called {\it the canonical symplectic
section} associated with the Lie algebroid $E$.

Suppose that $(x^i)$ are local coordinates on an open subset $U$
of $M$ and that $\{e_{\alpha}\}$ is a local basis of sections of
the vector bundle $\tau_{E}^{-1}(U)\rightarrow U$ as above. Then,
$\{\tilde{e}_{\alpha},\bar{e}_{\alpha}\}$ is a local basis of
sections of the vector bundle
$(\tau_{E}^{\tau_{E}^*})^{-1}((\tau_{E}^*)^{-1}(U))\rightarrow
(\tau_{E}^*)^{-1}(U)$, where $\tau_{E}^{\tau_{E}^*}:{\mathcal
T}^EE^*\rightarrow E^*$ is the vector bundle projection and
$$\tilde{e}_{\alpha}(a^*)=(e_{\alpha}(\tau_{E}^*(a^*)),\rho_{\alpha}^i\displaystyle\frac{\partial}{\partial
x^i}_{|a^*}),\;\;\;\bar{e}_{\alpha}(a^*)=(0,\displaystyle\frac{\partial}{\partial
y_{\alpha}}_{|a^*}).
$$
Here, $(x^i,y_{\alpha})$ are the local coordinates on $E^*$
induced by the local coordinates $(x^i)$ and the dual basis
$\{e^{\alpha}\}$ of $\{e_{\alpha}\}$. Moreover, we have that
$$
\begin{array}{c}
\lcf\tilde{e}_{\alpha},\tilde{e}_{\beta}\rcf_{E}^{\tau_{E}^*}\kern-2pt=C_{\alpha\beta}^{\gamma}\tilde{e}_{\gamma},\;\;
\lcf\tilde{e}_{\alpha},\bar{e}_{\beta}\rcf_{E}^{\tau_{E}^*}\kern-2pt=\lcf\bar{e}_{\alpha},\bar{e}_{\beta}\rcf_{E}^{\tau_{E}^*}\kern-2pt=0,\;\;
\rho_{E}^{\tau_{E}^*}(\tilde{e}_{\alpha})\kern-3pt=\rho_{\alpha}^i\displaystyle\frac{\partial}{\partial
x^i},\;\;\rho_{E}^{\tau_{E}^*}(\bar{e}_{\alpha})\kern-3pt=\displaystyle\frac{\partial}{\partial
y_{\alpha}},
\end{array}$$
and
\begin{equation}\label{formas}
\lambda_E(x^i,y_{\alpha})=y_{\alpha}\tilde{e}^{\alpha},\;\;\;\Omega_E(x^i,y_{\alpha})=\tilde{e}^{\alpha}\wedge\bar{e}^{\alpha}+\displaystyle\frac{1}{2}C_{\alpha\beta}^{\gamma}y_{\gamma}\tilde{e}^{\alpha}\wedge\tilde{e}^{\beta},
\end{equation}
(for more details, see \cite{LMM,M2}).

\subsection{Lie affgebroids}\label{sec1.2}

Let $\tau_A:A\rightarrow M$ be an affine bundle with associated
vector bundle $\tau_V:V\rightarrow M$. Denote by
$\tau_{A^+}:A^+=Aff(A,\R)\rightarrow M$ the dual bundle whose
fibre over $x\in M$ consists of affine functions on the  fibre
$A_x$. Note that this bundle has a distinguished section
$1_A\in\Gamma(\tau_{A^+})$ corresponding to the constant function
$1$ on $A$. We also consider the bidual bundle
$\tau_{\widetilde{A}}:\widetilde{A}\rightarrow M$ whose fibre at
$x\in M$ is the vector space $\widetilde{A}_x=(A_x^+)^*$. Then,
$A$ may be identified with an affine subbundle of $\widetilde{A}$
via the inclusion $i_A:A\rightarrow\widetilde{A}$ given by
$i_A(a)(\varphi)=\varphi(a)$, which is injective affine map whose
associated linear map is denoted by
$i_V:V\rightarrow\widetilde{A}$. Thus, $V$ may be identified with
a vector subbundle of $\widetilde{A}$. Using these facts, one can
prove that there is a one-to-one correspondence between affine
functions on $A$ and linear functions on $\widetilde{A}$. On the
other hand, there is an obvious one-to-one correspondence between
affine functions on $A$ and sections of $A^+$.

A {\it
 Lie affgebroid structure} on $A$ consists of a Lie algebra structure
 $\lcf\cdot,\cdot\rcf_V$ on the space
 $\Gamma(\tau_V)$ of the sections of $\tau_V:V\rightarrow M$, a $\R$-linear action
 $D:\Gamma(\tau_A)\times\Gamma(\tau_V)\rightarrow\Gamma(\tau_V)$ of
 the sections of $A$ on $\Gamma(\tau_V)$ and an affine map
 $\rho_A:A\rightarrow TM$, the {\it anchor map}, satisfying the following
 conditions:

\smallskip

  \begin{enumerate}
\item[$ \bullet$]
$D_X\lcf\bar{Y},\bar{Z}\rcf_V=\lcf
D_X\bar{Y},\bar{Z}\rcf_V+\lcf\bar{Y},D_X\bar{Z}\rcf_V,$
\item[$\bullet$]
$D_{X+\bar{Y}}\bar{Z}=D_X\bar{Z}+\lcf\bar{Y},\bar{Z}\rcf_V,$
\item[$\bullet$] $D_X(f\bar{Y})=fD_X\bar{Y}+\rho_A(X)(f)\bar{Y},$
\end{enumerate}

\smallskip

\noindent for $X\in\Gamma(\tau_A)$,
$\bar{Y},\bar{Z}\in\Gamma(\tau_V)$ and $f\in C^{\infty}(M)$ (see
\cite{GGrU,MMeS}).

If $(\lcf\cdot,\cdot\rcf_V,D,\rho_A)$ is a Lie affgebroid
structure on an affine bundle $A$ then
$(V,\lcf\cdot,\cdot\rcf_V,\rho_V)$ is a Lie algebroid, where
$\rho_V:V\rightarrow TM$ is the vector bundle map associated with
the affine morphism $\rho_A:A\rightarrow TM$.

A Lie affgebroid structure on an affine bundle
$\tau_A:A\rightarrow M$ induces a Lie algebroid structure
$(\lcf\cdot,\cdot\rcf_{\widetilde{A}},\rho_{\widetilde{A}})$ on
the bidual bundle $\widetilde{A}$ such that
$1_A\in\Gamma(\tau_{A^+})$ is a $1$-cocycle in the corresponding
Lie algebroid cohomology, that is, $d^{\widetilde{A}}1_A=0$.
Indeed, if $X_0\in\Gamma(\tau_A)$ then for every section
$\widetilde{X}$ of $\widetilde{A}$ there exists a unique function
$f\in C^{\infty}(M)$ and a unique section
$\bar{X}\in\Gamma(\tau_V)$ such that $\widetilde{X}=fX_0+\bar{X}$
and
$$
\begin{array}{rcl}
\rho_{\widetilde{A}}(fX_0+\bar{X})&=&f\rho_A(X_0)+\rho_V(\bar{X}),\\
\lcf
fX_0+\bar{X},gX_0+\bar{Y}\rcf_{\widetilde{A}}&=&(\rho_V(\bar{X})(g)-\rho_V(\bar{Y})(f)+f\rho_A(X_0)(g)\\&&-g\rho_A(X_0)(f))X_0
+\lcf\bar{X},\bar{Y}\rcf_V+fD_{X_0}\bar{Y}-gD_{X_0}\bar{X}.
\end{array}
$$

Conversely, let $(U,\lcf\cdot,\cdot\rcf_U,\rho_U)$ be a Lie
algebroid over $M$ and $\phi:U\rightarrow\R$ be a $1$-cocycle of
$(U,\lcf\cdot,\cdot\rcf_U,\rho_U)$ such that $\phi_{|U_x}\neq 0$,
for all $x\in M$. Then, $A=\phi^{-1}\{1\}$  is an affine bundle
over $M$ which admits a Lie affgebroid structure in such a way
that $(U,\lcf\cdot,\cdot\rcf_U,\rho_U)$ may be identified with the
bidual Lie algebroid
$(\widetilde{A},\lcf\cdot,\cdot\rcf_{\widetilde{A}},\rho_{\widetilde{A}})$
to $A$ and, under this identification, the $1$-cocycle
$1_A:\widetilde{A}\rightarrow\R$ is just $\phi$. The affine bundle
$\tau_A:A\rightarrow M$ is modelled on the vector bundle
$\tau_V:V=\phi^{-1}\{0\}\rightarrow M$. In fact, if
$i_V:V\rightarrow U$ and $i_A:A\rightarrow U$ are the canonical
inclusions, then
$$
\begin{array}{rcl} i_V\circ\lcf\bar{X},\bar{Y}\rcf_V&=&\lcf
i_V\circ\bar{X},i_V\circ\bar{Y}\rcf_U,\;\;i_V\circ D_X\bar{Y}=\lcf
i_A\circ X,i_V\circ\bar{Y}\rcf_U,\\
\rho_A(X)&=&\rho_{U}(i_A\circ X),
\end{array}
$$
 for $\bar{X},\bar{Y}\in\Gamma(\tau_V)$ and $X\in\Gamma(\tau_A)$
 (for more details, see \cite{GGrU,MMeS}).

A trivial example of a Lie affgebroid may be constructed as
follows. Let $\tau:M\to\R$ be a fibration and
$\tau_{1,0}:J^1\tau\to M$ be the $1$-jet bundle of local sections
of $\tau:M\to\R$. It is well known that $\tau_{1,0}:J^1\tau\to M$
is an affine bundle modelled on the vector bundle
$\pi=(\pi_M)_{|V\tau}:V\tau\to M$, where $V\tau$ is the vertical
bundle of $\tau:M\to\R$. Moreover, if $t$ is the usual coordinate
on $\R$ and $\eta$ is the closed $1$-form on $M$ given by
$\eta=\tau^*(dt)$ then we have the following identification
$$J^1\tau\cong\{v\in TM/\eta(v)=1\}$$
(see, for instance, \cite{Sa}). Note that $V\tau=\{v\in
TM/\eta(v)=0\}.$ Thus, the bidual bundle $\widetilde{J^1\tau}$ to
the affine bundle $\tau_{1,0}:J^1\tau\to M$ may be identified with
the tangent bundle $TM$ to $M$ and, under this identification, the
Lie algebroid structure on $\pi_M:TM\to M$ is the standard Lie
algebroid structure and the $1$-cocycle $1_{J^1\tau}$ on
$\pi_M:TM\to M$ is just the closed $1$-form $\eta$.

\setcounter{equation}{0}
\section{Hamiltonian formalism on Lie affgebroids}

\subsection{A cosymplectic structure on ${\mathcal
T}^{\widetilde{A}}V^*$}

Suppose that $(\tau_A:A\rightarrow M, \tau_V:V\rightarrow M,
(\lcf\cdot,\cdot\rcf_V,$ $D,\rho_A))$ is a Lie affgebroid and let
$({\mathcal
T}^{\widetilde{A}}V^*,\lcf\cdot,\cdot\rcf_{\widetilde{A}}^{\tau_{V}^*},\rho_{\widetilde{A}}^{\tau_{V}^*})$
be the prolongation of the bidual Lie algebroid
$(\widetilde{A},\lcf\cdot,\cdot\rcf_{\widetilde{A}},$ $
\rho_{\widetilde{A}})$ over the fibration
$\tau_{V}^*:V^*\rightarrow M$. In this section, we are going to
construct a cosymplectic structure on ${\mathcal
T}^{\widetilde{A}}V^*$ using the canonical symplectic section
associated with the Lie algebroid $\widetilde{A}$ and a
Hamiltonian section (this construction was done in \cite{IMPS}).

Let $(x^i)$ be local coordinates on an open subset $U$ of $M$ and
$\{e_0,e_{\alpha}\}$ be a local basis of sections of the vector
bundle $\tau^{-1}_{\widetilde{A}}(U)\rightarrow U$ adapted to the
$1$-cocycle $1_A$ (that is, $1_A(e_0)=1$ and $1_A(e_\alpha)=0$,
for all $\alpha$) and such that
\begin{equation}\label{coranc}
\begin{array}{rclrclccrclrcl}
\lcf
e_0,e_{\alpha}\rcf_{\widetilde{A}}=C_{0\alpha}^{\gamma}e_{\gamma},\;\;\lcf
e_{\alpha},e_{\beta}\rcf_{\widetilde{A}}=C_{\alpha\beta}^{\gamma}e_{\gamma},\;\;
\rho_{\widetilde{A}}(e_0)=\rho_0^i\displaystyle\frac{\partial}{\partial
x^i},\;\;
\rho_{\widetilde{A}}(e_{\alpha})=\rho_{\alpha}^i\displaystyle\frac{\partial}{\partial
x^i}.
\end{array}
\end{equation}

Denote by $(x^i,y^0,y^{\alpha})$ the corresponding local
coordinates on $\widetilde{A}$ and by $(x^i,y_0,y_{\alpha})$ the
dual coordinates on the dual vector bundle $\tau_{A^+}:A^+\to M$
to $\widetilde{A}$. Then, $(x^i,y_{\alpha})$ are local coordinates
on $V^*$ and $\{\tilde{e}_0,\tilde{e}_{\alpha},\bar{e}_{\alpha}\}$
is a local basis of sections of the vector bundle
$\tau^{\tau_{V}^*}_{\widetilde{A}}:{\mathcal
T}^{\widetilde{A}}V^*\rightarrow V^*$, where
\begin{equation}\label{tilbar}
\tilde{e}_0(\psi)=(e_0(\tau_V^*(\psi)),\rho_0^i\displaystyle\frac{\partial}{\partial
x^i}_{|\psi}),\;\;\tilde{e}_{\alpha}(\psi)=(e_{\alpha}(\tau_V^*(\psi)),\rho_{\alpha}^i\displaystyle\frac{\partial}{\partial
x^i}_{|\psi}),\;\;\bar{e}_{\alpha}(\psi)=(0,\displaystyle\frac{\partial}{\partial
y_{\alpha}}_{|\psi}). \end{equation} Using this local basis one
may introduce local coordinates $(x^i,y_{\alpha};z^0,$
$z^{\alpha},v_{\alpha})$ on ${\mathcal T}^{\widetilde{A}}V^*$. A
direct computation proves that
$$
\begin{array}{l}
\rho_{\widetilde{A}}^{\tau_V^*}(\tilde{e}_0)=\rho_0^i\displaystyle\frac{\partial
}{\partial
x^i},\;\;\;\rho_{\widetilde{A}}^{\tau_V^*}(\tilde{e}_\alpha)=\rho_\alpha^i\displaystyle\frac{\partial
}{\partial x^i},\;\;\;
\rho_{\widetilde{A}}^{\tau_V^*}(\bar{e}_\alpha)=\displaystyle\frac{\partial
}{\partial y_\alpha},\\
\lcf\tilde{e}_0,\tilde{e}_\beta\rcf_{\widetilde{A}}^{\tau_V^*}=C_{0\beta}^\gamma\tilde{e}_\gamma,\;\;\;
\lcf\tilde{e}_\alpha,\tilde{e}_\beta\rcf_{\widetilde{A}}^{\tau_V^*}=C_{\alpha\beta}^\gamma\tilde{e}_\gamma,\\
\lcf\tilde{e}_0,\bar{e}_\alpha\rcf_{\widetilde{A}}^{\tau_V^*}=\lcf\tilde{e}_\alpha,\bar{e}_\beta\rcf_{\widetilde{A}}^{\tau_V^*}=\lcf
\bar{e}_\alpha, \bar{e}_\beta\rcf_{\widetilde{A}}^{\tau_V^*}=0,
\end{array}
$$
for all $\alpha$ and $\beta$. Thus, if
$\{\tilde{e}^0,\tilde{e}^\alpha,\bar{e}^\alpha\}$ is the dual
basis of $\{\tilde{e}_0,\tilde{e}_\alpha,\bar{e}_\alpha\}$ then
\begin{equation}\label{dif*}
\begin{array}{rcl}
d^{{\mathcal T}^{\widetilde{A}}V^*}f&=&
\rho_0^i\displaystyle\frac{\partial f}{\partial x^i}\tilde{e}^0+
\rho_\alpha^i\displaystyle\frac{\partial f}{\partial
x^i}\tilde{e}^\alpha + \displaystyle\frac{\partial f}{\partial
y_\alpha} \bar{e}^\alpha,\\
d^{{\mathcal
T}^{\widetilde{A}}V^*}\tilde{e}^\gamma&=&\displaystyle
-\frac{1}{2}C_{0\alpha}^\gamma\tilde{e}^0\wedge\tilde{e}^\alpha
\displaystyle -\frac{1}{2} C_{\alpha\beta}^{\gamma}
\tilde{e}^\alpha \wedge \tilde{e}^\beta, \makebox[1cm]{}
d^{{\mathcal T}^{\widetilde{A}}V^*}\tilde{e}^0=d^{{\mathcal
T}^{\widetilde{A}}V^*}\bar{e}^\gamma=0,
\end{array}
\end{equation}
for $f\in C^\infty(V^*).$

Let $\mu:A^+\rightarrow V^*$ be the canonical projection given by
$\mu(\varphi)=\varphi^l$, for $\varphi\in A^+_x$, with $x\in M$,
where $\varphi^l\in V^*_x$ is the linear map associated with the
affine map $\varphi$ and $h:V^*\rightarrow A^+$ be a Hamiltonian
section of $\mu$, that is, $\mu\circ h=Id$.

Now, we consider the Lie algebroid prolongation ${\mathcal
T}^{\widetilde{A}}A^+$ of the Lie algebroid $\widetilde{A}$ over
$\tau_{A^+}:A^+\to M$ with vector bundle projection
$\tau^{\tau_{A^+}}_{\widetilde{A}}:{\mathcal
T}^{\widetilde{A}}A^+\rightarrow A^+$ (see Section
\ref{sec1.1.1}). Then, we may introduce the map ${\mathcal
T}h:{\mathcal T}^{\widetilde{A}}V^*\rightarrow{\mathcal
T}^{\widetilde{A}}A^+$ defined by ${\mathcal
T}h(\tilde{a},X_{\alpha})=(\tilde{a},(T_{\alpha}h)(X_{\alpha})),$
for $(\tilde{a},X_{\alpha})\in({\mathcal
T}^{\widetilde{A}}V^*)_{\alpha}$, with $\alpha\in V^*.$ It is easy
to prove that the pair $({\mathcal T}h,h)$ is a Lie algebroid
morphism between the Lie algebroids
$\tau_{\widetilde{A}}^{\tau_V^*}:{\mathcal
T}^{\widetilde{A}}V^*\rightarrow V^*$ and
$\tau_{\widetilde{A}}^{\tau_{A^+}}:{\mathcal
T}^{\widetilde{A}}A^+\rightarrow A^+$.

\noindent Next, denote by $\lambda_h$ and $\Omega_h$ the sections
of the vector bundles $({\mathcal
T}^{\widetilde{A}}V^*)^*\rightarrow V^*$ and $\Lambda^2({\mathcal
T}^{\widetilde{A}}V^*)^*\rightarrow V^*$ given by
\begin{equation}\label{Omegah}
\lambda_h=({\mathcal
T}h,h)^*(\lambda_{\widetilde{A}}),\;\;\Omega_h=({\mathcal
T}h,h)^*(\Omega_{\widetilde{A}}), \end{equation}

where $\lambda_{\widetilde{A}}$ and $\Omega_{\tilde{A}}$ are the
Liouville section and  the canonical symplectic section,
respectively, associated with the Lie algebroid $\widetilde{A}.$
Note that $\Omega_h=-d^{{\mathcal
T}^{\widetilde{A}}V^*}\lambda_h.$

On the other hand, let $\eta:{\mathcal
T}^{\widetilde{A}}V^*\rightarrow\R$ be the section of $({\mathcal
T}^{\widetilde{A}}V^*)^*\rightarrow V^*$ defined by
\begin{equation}\label{eta}
\eta(\tilde{a},X_{\alpha})=1_A(\tilde{a}),
\end{equation}
 for
$(\tilde{a},X_{\alpha})\in({\mathcal
T}^{\widetilde{A}}V^*)_{\alpha}$, with $\alpha\in V^*$. Note that
if $pr_1:{\mathcal T}^{\widetilde{A}}V^*\to \widetilde{A}$ is the
canonical projection on the first factor then $(pr_1,\tau_V^*)$ is
a morphism between the Lie algebroids
$\tau_{\widetilde{A}}^{\tau_V^*}:{\mathcal
T}^{\widetilde{A}}V^*\rightarrow V^*$ and
$\tau_{\widetilde{A}}:\widetilde{A}\to M$ and
$(pr_1,\tau_V^*)^*(1_A)=\eta$. Thus, since $1_A$ is a $1$-cocycle
of $\tau_{\widetilde{A}}:\widetilde{A}\rightarrow M$, we deduce
that $\eta$ is a $1$-cocycle of the Lie algebroid
$\tau_{\widetilde{A}}^{\tau_V^*}:{\mathcal
T}^{\widetilde{A}}V^*\rightarrow V^*.$

Suppose that
$h(x^i,y_{\alpha})=(x^i,-H(x^j,y_{\beta}),y_{\alpha})$. Then
$\eta=\tilde{e}^0$ and, from (\ref{formas}), (\ref{Omegah}) and
the definition of the map ${\mathcal T}h$, it follows that
\begin{equation}\label{Omegahlocal}
\Omega_h=\tilde{e}^{\gamma}\wedge\bar{e}^{\gamma}+\frac{1}{2}C_{\gamma\beta}^{\alpha}y_{\alpha}\tilde{e}^{\gamma}\wedge\tilde{e}^{\beta}+(\rho_{\gamma}^i\frac{\partial
H}{\partial
x^i}-C_{0\gamma}^{\alpha}y_{\alpha})\tilde{e}^{\gamma}\wedge\tilde{e}^0+\frac{\partial
H }{\partial y_{\gamma}}\bar{e}^{\gamma}\wedge\tilde{e}^0.
\end{equation}

Thus, it is easy to prove that the pair $(\Omega_h,\eta)$ is a
cosymplectic structure on the Lie algebroid
$\tau_{\widetilde{A}}^{\tau_V^*}:{\mathcal
T}^{\widetilde{A}}V^*\rightarrow V^*$, that is,
$$
\begin{array}{c}
\{\eta\wedge\Omega_h\wedge\dots^{(n}
\dots\wedge\Omega_h\}(\alpha)\neq 0,\makebox[1.5cm]{for
all}\alpha\in V^*,\\
d^{{\mathcal T}^{\widetilde{A}}V^*}\eta=0,\;\;\;d^{{\mathcal
T}^{\widetilde{A}}V^*}\Omega_h=0.
\end{array}
$$

\begin{remark} {\em Let ${\mathcal T}^VV^*$ be the prolongation of
the Lie algebroid $V$ over the projection $\tau_V^*:V^*\to M$.
Denote by $\lambda_V$ and $\Omega_V$ the Liouville section and the
canonical symplectic section, respectively, of $V$ and by
$(i_V,Id):{\mathcal T}^VV^*\to{\mathcal T}^{\widetilde{A}}V^*$ the
canonical inclusion. Then, using (\ref{lambdaE}), (\ref{Omegah}),
(\ref{eta}) and the fact that $\mu\circ h=Id$, we obtain that
$$(i_V,Id)^*(\lambda_h)=\lambda_V,\;\;\;(i_V,Id)^*(\eta)=0.$$
Thus, since $(i_V,Id)$ is a Lie algebroid morphism over the
identity of $V^*$, we also deduce that
$$(i_V,Id)^*(\Omega_h)=\Omega_V.$$
\begin{flushright}$\diamondsuit$\end{flushright}}
\end{remark}

Now, given a section $\gamma$ of the dual bundle $V^*$ to $V$, we
can consider the morphism $({\mathcal T}\gamma,\gamma)$ between
the vector bundles $\tau_{\widetilde{A}}:\widetilde{A}\to M$ and
$\tau_{\widetilde{A}}^{\tau_V^*}:{\mathcal
T}^{\widetilde{A}}V^*\to V^*$

\begin{picture}(375,90)(40,10)
\put(190,20){\makebox(0,0){$M$}}
\put(250,25){$\gamma$}\put(210,20){\vector(1,0){80}}
\put(310,20){\makebox(0,0){$V^*$}}
\put(170,50){$\tau_{\widetilde{A}}$}
\put(190,70){\vector(0,-1){40}}
\put(320,50){$\tau_{\widetilde{A}}^{\tau_V^*}$}
\put(310,70){\vector(0,-1){40}}
\put(190,80){\makebox(0,0){$\widetilde{A}$}}
\put(245,85){${\mathcal T}\gamma$}\put(210,80){\vector(1,0){80}}
\put(310,80){\makebox(0,0){${\mathcal T}^{\widetilde{A}}V^*$}}
\end{picture}

defined by ${\mathcal
T}\gamma(\tilde{a})=(\tilde{a},(T_x\gamma)(\rho_{\widetilde{A}}(\tilde{a}))),$
for $\tilde{a}\in\widetilde{A}_x$ and $x\in M$.

\begin{theorem}\label{thmor}If $\gamma$ is a section of the vector bundle $\tau_V^*:V^*\to M$ then the pair
$({\mathcal T}\gamma,\gamma)$ is a morphism between the Lie
algebroids
$(\widetilde{A},\lcf\cdot,\cdot\rcf_{\widetilde{A}},\rho_{\widetilde{A}})$
and $({\mathcal T}^{\widetilde{A}}V^*,
\lcf\cdot,\cdot\rcf_{\widetilde{A}}^{\tau_V^*},\rho_{\widetilde{A}}^{\tau_V^*}).$
Moreover,
$$({\mathcal T}\gamma,\gamma)^*\lambda_h=h\circ\gamma,\;\;\;
({\mathcal
T}\gamma,\gamma)^*(\Omega_h)=-d^{\widetilde{A}}(h\circ\gamma).$$
\end{theorem}
\begin{proof} Suppose that $(x^i)$ are local coordinates on $M$, that
$\{e_0,e_\alpha\}$ is a local basis of
$\Gamma(\tau_{\widetilde{A}})$ adapted to $1_A$ and that
\[
\gamma=\gamma_\alpha e^\alpha,
\]
with $\gamma_\alpha$ local real functions on $M$ and
$\{e^0,e^\alpha\}$ the dual basis to $\{e_0,e_\alpha\}$. Denote by
$\{\tilde{e}_0,\tilde{e}_\alpha,\bar{e}_\alpha\}$ the
corresponding local basis of
$\Gamma(\tau_{\widetilde{A}}^{\tau_V^*})$. Then, using
(\ref{tilbar}), it follows that
$${\mathcal T}\gamma\circ
e_0=(\tilde{e}_0+\rho_0^{i}\frac{\partial\gamma_{\nu}}{\partial
x^i}\bar{e}_{\nu})\circ\gamma,\;\;{\mathcal T}\gamma\circ
e_{\alpha}=(\tilde{e}_\alpha + \rho_\alpha^i\frac{\partial
\gamma_\nu}{\partial x^i}\bar{e}_\nu)\circ \gamma ,$$
for $\alpha\in \{1,\dots ,n\}$, $\rho_0^i,\rho_\alpha^i$ being the
components of the anchor map of $\widetilde{A}$ with respect  to
the local coordinates $(x^i)$ and to the basis $\{e_0,e_\alpha\}.$
Thus,
\[
\begin{array}{rcl}
({\mathcal T}\gamma,\gamma)^*(\tilde{e}^0)&=&e^0,\;\;({\mathcal T}\gamma,\gamma)^*(\tilde{e}^{\alpha})=e^{\alpha},\\
({\mathcal
T}\gamma,\gamma)^*(\bar{e}^\alpha)&=&\rho_0^i\displaystyle\frac{\partial\gamma_{\alpha}}{\partial
x^i}e^0+ \rho_\beta^i\displaystyle\frac{\partial
\gamma_\alpha}{\partial x^i}e^\beta = d^{\widetilde{A}}
\gamma_{\alpha},
\end{array}
\]
where $\{\tilde{e}^0,\tilde{e}^\alpha,\bar{e}^\alpha\}$ is the
dual basis to $\{\tilde{e}_0,\tilde{e}_\alpha,\bar{e}_\alpha\}.$

Therefore, from (\ref{dE}), (\ref{coranc}) and (\ref{dif*}), we
obtain that the pair $({\mathcal T}\gamma,\gamma)$ is a morphism
between the Lie algebroids $\widetilde{A}\to M$ and ${\mathcal
T}^{\widetilde{A}}V^*\to V^*.$

On the other hand, if $x$ is a point of $M$ and
$\tilde{a}\in\widetilde{A}_x$ then, using (\ref{lambdaE}), we have
$$
\begin{array}{rcl}
(({\mathcal
T}\gamma,\gamma)^*\lambda_h)(x)(\tilde{a})&=&(({\mathcal
T}\gamma,\gamma)^*({\mathcal
T}h,h)^*\lambda_{\widetilde{A}})(x)(\tilde{a})=(({\mathcal
T}(h\circ\gamma),h\circ\gamma)^*\lambda_{\widetilde{A}})(x)(\tilde{a})\\
&=&\lambda_{\widetilde{A}}(h(\gamma(x)))(\tilde{a},T(h\circ\gamma)(\rho_{\widetilde{A}}(\tilde{a})))=(h\circ\gamma)(x)(\tilde{a})
\end{array}
$$
that is, $$({\mathcal T}\gamma,\gamma)^*\lambda_h=h\circ\gamma.$$
Consequently, since $\Omega_h=-d^{{\mathcal
T}^{\widetilde{A}}V^*}\lambda_h$ and $({\mathcal T}\gamma,\gamma)$
is a morphism between the Lie algebroids $\widetilde{A}$ and
${\mathcal T}^{\widetilde{A}}V^*$, we deduce that
$$({\mathcal
T}\gamma,\gamma)^*(\Omega_h)=-d^{\widetilde{A}}(h\circ\gamma).$$

\end{proof}

\subsection{The Hamilton equations on a Lie affgebroid}

Let $h:V^*\to A^+$ be a Hamiltonian section of $\mu:A^+\to V^*$
and $R_h\in\Gamma(\tau_{\widetilde{A}}^{\tau_V^*})$ be the Reeb
section of the cosymplectic structure $(\Omega_h,\eta)$ on
${\mathcal T}^{\widetilde{A}}V^*$. $R_h$ is characterized by the
following conditions
$$i_{R_h}\Omega_h=0\makebox[1cm]{and}i_{R_h}\eta=1.$$
With respect to the basis
$\{\tilde{e}_0,\tilde{e}_{\alpha},\bar{e}_{\alpha}\}$ of
$\Gamma(\tau_{\widetilde{A}}^{\tau_V^*})$, $R_h$ is locally
expressed as follows:
$$
R_h=\tilde{e}_0+\frac{\partial H}{\partial
y_{\alpha}}\tilde{e}_{\alpha}-(C_{\alpha\beta}^{\gamma}y_{\gamma}\frac{\partial
H}{\partial y_{\beta}}+\rho^i_{\alpha}\frac{\partial H}{\partial
x^i}-C_{0\alpha}^{\gamma}y_{\gamma})\bar{e}_{\alpha}.
$$

Thus, the integral sections of $R_h$ (i.e., the integral curves of
the vector field $\rho_{\widetilde{A}}^{\tau_{V}^*}(R_h)$) satisfy
the following equations
$$
\frac{dx^i}{dt}=\rho_0^i+\frac{\partial H}{\partial
y_{\alpha}}\rho_{\alpha}^i,\;\;\;\frac{dy_{\alpha}}{dt}=-\rho_{\alpha}^i\frac{\partial
H }{\partial
x^i}+y_{\gamma}(C_{0\alpha}^{\gamma}+C_{\beta\alpha}^{\gamma}\frac{\partial
H }{\partial y_{\beta}}),
$$
for $i\in\{1,\dots,m\}$ and $\alpha\in\{1,\dots,n\}$. These
equations are called the {\it Hamilton equations for $h$} (see
\cite{M}) and the section $R_h$ is called {\it the Hamiltonian
section associated with $h$} (see \cite{IMPS}).

\setcounter{equation}{0}
\section{The Hamilton-Jacobi equation on Lie affgebroids}
In this section, we will prove the main result of the paper.

Let $\tau_A:A\to M$ be an affine bundle with associated vector
bundle $\tau_V:V\to M$ and suppose that $h:V^*\to A^+$ is a
section of the canonical projection $\mu:A^+\to V^*$ and that
$\alpha:M\to A^+$ is a section of the vector bundle
$\tau_{A^+}:A^+\to M$. Then, $h\circ\mu\circ\alpha:M\to A^+$ is
also a section of the vector bundle $\tau_{A^+}:A^+\to M$ and
there exists a unique real function on $M$, which we will denote
by $f(h,\alpha)$, such that
\begin{equation}\label{*}
\alpha-h\circ\mu\circ\alpha=f(h,\alpha)1_A.
\end{equation}
Now, we will prove the result announced above.

\begin{theorem}\label{thHJ} Let $\tau_A:A\rightarrow M$ be a Lie affgebroid modelled over the vector bundle $\tau_V:V\rightarrow M$
with Lie affgebroid structure $(\lcf\cdot,\cdot\rcf_V,D,\rho_A)$
and $({\mathcal
T}^{\widetilde{A}}V^*,\lcf\cdot,\cdot\rcf_{\widetilde{A}}^{\tau_{V}^*},\rho_{\widetilde{A}}^{\tau_{V}^*})$
be the prolongation of the bidual Lie algebroid
$(\widetilde{A},\lcf\cdot,\cdot\rcf_{\widetilde{A}},\rho_{\widetilde{A}})$
over the fibration $\tau_{V}^*:V^*\rightarrow M$. Suppose that
$h:V^*\to A^+$ is a Hamiltonian section of the canonical
projection $\mu:A^+\to V^*$ and that
$R_h\in\Gamma(\tau_{\widetilde{A}}^{\tau_V^*})$ is the
corresponding Hamiltonian section. Let
$\alpha\in\Gamma(\tau_{A^+})$ be a $1$-cocycle,
$d^{\widetilde{A}}\alpha=0$, and denote by
$R_h^{\alpha}\in\Gamma(\tau_{\widetilde{A}})$ the section of
$\tau_{\widetilde{A}}:\widetilde{A}\to M$ given by
$R_h^{\alpha}=pr_1\circ R_h\circ \mu\circ\alpha$, where
$pr_1:{\mathcal T}^{\widetilde{A}}V^*\to\widetilde{A}$ is the
canonical projection over the first factor. Then, the two
following conditions are equivalent:
\begin{enumerate}
\item For every integral curve $c$ of the vector field $\rho_{\widetilde{A}}(R_h^\alpha)$, that is,
$c$ is a curve on $M$ such that
\begin{equation}\label{eqc}
\rho_{\widetilde{A}}(R_h^{\alpha})(c(t))=\dot{c}(t),\makebox[1.2cm]{for
all}t,
\end{equation}
the curve $t\to (\mu\circ\alpha\circ c)(t)$ on $V^*$ satisfies the
Hamilton equations for $h$.
\item $\alpha$ satisfies the Hamilton-Jacobi equation
$d^V(f(h,\alpha))=0$, that is, $f(h,\alpha)$ is constant on the
leaves of the Lie algebroid foliation of $V$.
\end{enumerate}
\end{theorem}

\begin{proof}
For a curve $c:I=(-\epsilon,\epsilon)\subset\R\to M$ on the base
we define the curves $\map{\beta}{I}{V^*}$ and
$\map{\gamma}{I}{\widetilde{A}}$ by
$$
\beta(t)=\mu(\alpha(c(t)))\makebox[1.5cm]{and}\gamma(t)=R_h^{\alpha}(c(t)).
$$
Since $\mu\circ\alpha$ and $R_h^{\alpha}$ are sections of
$\tau_{V}^*: V^* \to M$ and $\tau_{\widetilde{A}}: \widetilde{A}
\to M$, respectively, it follows that both curves project to $c$.

We consider the curve $v=(\gamma,\dot{\beta})$ in
$\widetilde{A}\times TV^*$ and notice the following important
facts about $v$:
\begin{itemize}
\item $v(t)$ is in ${\mathcal T}^{\widetilde{A}}V^*$, for every $t\in I$, if and
only if $c$ satisfies (\ref{eqc}). Indeed
$\rho_{\widetilde{A}}\circ\gamma=\rho_{\widetilde{A}}\circ
R_h^{\alpha}\circ c$ while $T\tau_V^*\circ\dot{\beta}=\dot{c}$.
\item In such a case, $\beta$ is a solution of the Hamilton equations for $h$ if and only if $v(t)=R_h(\beta(t))$,
for every $t\in I$. Indeed, the first components coincide
$pr_1(v(t))=\gamma(t)$ and
$pr_1(R_h(\beta(t)))=pr_1(R_h(\mu(\alpha(c(t)))))
=R_h^{\alpha}(c(t))=\gamma(t)$, and the equality of the second
components is just
$\dot{\beta}(t)=\rho_{\widetilde{A}}^{\tau_V^*}(R_h(\beta(t)))$.
\end{itemize}

We denote by $\alpha_\mu=\mu\circ\alpha$ which is a section of
$\tau_V^*:V^*\to M$ and we also consider the map ${\mathcal
T}\alpha_\mu:\widetilde{A}\to {\mathcal T}^{\widetilde{A}}V^*$
given by ${\mathcal
T}\alpha_\mu=(Id,T(\alpha_{\mu})\circ\rho_{\widetilde{A}})$. We
recall that $({\mathcal
T}\alpha_\mu,\alpha_\mu)^*\Omega_h=-d^{\widetilde{A}}(h\circ\alpha_\mu)$,
because of Theorem \ref{thmor}.

\smallskip\noindent[$(ii)$ $\Rightarrow$ $(i)$]
Assume that $c$ satisfies (\ref{eqc}), so that $v(t)$ is a curve
in ${\mathcal T}^{\widetilde{A}}V^*$. We have to prove that $v(t)$
equals to $R_h(\beta(t))$, for every $t\in I$.

The difference $d(t)=v(t)-R_h(\beta(t))$ is vertical with respect
to the projection $pr_1:{\mathcal
T}^{\widetilde{A}}V^*\to\widetilde{A}$, that is, $pr_{1}(d(t)) =
0$, for all $t$ (note that $pr_{1}(v(t)) = pr_{1}(R_h(\beta(t))) =
\gamma (t)$, for all $t$). Therefore, we have that $\eta(d(t))=0$
and $\Omega_h (\beta(t))$ $(d(t), \varsigma(t))= 0$, for every
vertical curve $t \to \varsigma(t)$ (see (\ref{tilbar}) and
(\ref{Omegahlocal})).

Let $\map{\tilde{a}}{I}{\widetilde{A}}$ be any curve on
$\widetilde{A}$ over $c$ (that is, $\tau_{\widetilde{A}}
\circ\tilde{a} = c$). We consider its image under ${\mathcal
T}\alpha_\mu$, i.e., $\zeta(t)=({\mathcal
T}\alpha_{\mu})(\tilde{a}(t))=(\tilde{a}(t),T(\mu\circ\alpha)(\rho_{\widetilde{A}}(\tilde{a}(t))))$.
>From (\ref{eqc}), it follows that $v(t)=({\mathcal
T}\alpha_{\mu})(\gamma(t))$ is also in the image of ${\mathcal
T}\alpha_{\mu}$. Thus, using (\ref{*}) and the fact that $R_h$ is
the Reeb section of $(\Omega_h,\eta)$, we have that
$$
\begin{array}{rcl} \Omega_h (\beta(t))(d(t),\zeta(t)) \kern-5pt&=&\kern-5pt\Omega_h
(\beta(t))({\mathcal T}\alpha_{\mu}(\gamma(t)),{\mathcal
T}\alpha_{\mu}(\tilde{a}(t)))=-d^{\widetilde{A}}(h\circ\mu\circ\alpha)(c(t))(\gamma(t),\tilde{a}(t)) \\
\kern-5pt&=&\kern-5pt(d^{\widetilde{A}}(f(h,\alpha))\wedge
1_A)(c(t))(\gamma(t),\tilde{a}(t)).
\end{array}
$$
Now, since $d^V(f(h,\alpha))=0$ and the inclusion
$i_V:V\to\widetilde{A}$ is a Lie algebroid morphism, we deduce
that $d^{\widetilde{A}}(f(h,\alpha))\wedge 1_A=0$. This implies
that
$$\Omega_h(\beta(t))(d(t),\zeta(t))=0.$$

Since any element in $({\mathcal
T}^{\widetilde{A}}V^*)_{\alpha_\mu(x)}$, with $x \in M$, can be
obtained as a sum of an element in the image of ${\mathcal
T}\alpha_{\mu}$ and a vertical, we conclude that $\Omega_h
(\beta(t))(d(t),\nu(t))=\Omega_h (\mu(\alpha(c(t))))(d(t),\nu(t))
= 0$ for every curve $t \to \nu(t)$. Thus, since $\eta(d(t))=0$,
we deduce that $d(t)=0$.

\smallskip\noindent[$(i)$ $\Rightarrow$ $(ii)$]
Suppose that $x$ is a point of $M$ and that $b\in V_x$. We will
show that
$$\langle d^V(f(h,\alpha)),b\rangle=0.$$
Let $c: I = (-\epsilon, \epsilon) \to M$ be the integral curve of
$\rho_{\widetilde{A}}(R_h^{\alpha})$ such that $c(0)=x$. It
follows that $c$ satisfies (\ref{eqc}). Let
$\beta=\mu\circ\alpha\circ\ c$, $\gamma=R_h^{\alpha}\circ c$ and
$v=(\gamma,\dot{\beta})$ as above. Since $\beta$ satisfies the
Hamilton equations for $h$, we have that $v(t)=R_h(\beta(t))$ for
all $t$. Now, we take any curve $t \to a(t)$ in $V$ over $c$ such
that $a(0)=b$ and denote by $t\to\tilde{a}(t)$ its inclusion in
$\widetilde{A}$. Since $v(t)={\mathcal T}\alpha_{\mu}(\gamma(t))$
we have that
$$
\begin{array}{rclcl}
0\kern-5pt&=&\kern-5pt\Omega_h(\beta(t))(R_h(\beta(t)),{\mathcal
T}\alpha_{\mu}(\tilde{a}(t)))
\kern-5pt&=&\kern-5pt\Omega_h(\beta(t))(v(t),{\mathcal T}\alpha_{\mu}(\tilde{a}(t)))\\
\kern-5pt&=&\kern-5pt\Omega_h(\beta(t))({\mathcal
T}\alpha_{\mu}(\gamma(t)),{\mathcal T}\alpha_{\mu}(\tilde{a}(t)))
\kern-5pt&=&\kern-5pt-d^{\widetilde{A}}(h\circ\mu\circ\alpha)(c(t))(\gamma(t),\tilde{a}(t))\\
\kern-5pt&=&\kern-5pt(d^{\widetilde{A}}f(h,\alpha)\wedge
1_A)(c(t))(\gamma(t),\tilde{a}(t)).
\end{array}
$$
Thus, using that $1_A(c(t))(\gamma(t))=1$ (note that
$\eta(R_h)=1$) and that $1_A(\tilde{a}(t))=0$, we deduce that
$$0=(d^V(f(h,\alpha)))(c(t))(a(t)).$$
In particular, at $t=0$, we have that $\langle
d^V(f(h,\alpha)),b\rangle=0$.
\end{proof}

\begin{remark}
{\rm Obviously, we can consider as a cocycle $\alpha$ a
$1$-coboundary $\alpha=d^{\widetilde{A}}S$, for some function $S$
on $M$. Nevertheless, it should be noticed that on a Lie algebroid
there exist, in general, 1-cocycles that are not locally
$1$-coboundaries.}
\end{remark}

\setcounter{equation}{0}
\section{Examples}
\subsection{The Hamilton-Jacobi equation on Lie algebroids}\label{ex5.0}
Let $\tau_{E}: E \to M$ be a Lie algebroid over a manifold $M$.
Then, $\tau_{E}: E \to M$ is an affine bundle and the Lie
algebroid structure induces a Lie affgebroid structure on
$\tau_{E}: E \to M$. In fact, the dual bundle to $E$, as an affine
bundle, is the vector bundle $\tau_{E^+}: E^+ = E^* \times \R \to
M$, the bidual bundle is the vector bundle $\tau_{\widetilde{E}}:
\widetilde{E} = E \times \R \to M$ and the Lie algebroid structure
$(\lcf\cdot ,\cdot \rcf_{\widetilde{E}}, \rho_{\widetilde{E}})$ on
$\tau_{\widetilde{E}}: \widetilde{E} = E \times \R \to M$ is given
by
\[
\lcf (X, f), (Y, g) \rcf_{\widetilde{E}} = (\lcf X, Y\rcf_{E},
\rho_{E}(X)(g) - \rho_{E}(Y)(f)), \; \; \rho_{\widetilde{E}}(X, f)
= \rho_{E}(X),
\]
for $(X, f), (Y, g) \in \Gamma(\tau_{\widetilde{E}}) \cong
\Gamma(\tau_{E}) \times C^{\infty}(M)$, where $(\lcf\cdot ,
\cdot\rcf_{E}, \rho_{E})$ is the Lie algebroid structure on $E$.
The $1$-cocycle $1_{E}$ on $\widetilde{E}$ is the section $(0, 1)$
of $\tau_{E^+}: E^+ = E^* \times \R \to M$.

The map $\mu: E^+ = E^* \times \R \to E^*$ is the canonical
projection over the first factor and a Hamiltonian section $h: E^*
\to E^+ = E^* \times \R$ may be identified with a Hamiltonian
function $H$ on $E^*$ in such a way that
\[
h(\beta_{x}) = (\beta_{x}, -H(\beta_{x})), \; \; \mbox{ for }
\beta_{x} \in E^*_x \mbox{ and } x \in M.
\]
Now, if $\alpha$ is a $1$-cocycle of $\tau_{E}: E \to M$ then
$\alpha$ may be considered as the section $(\alpha, 0)$ of
$\tau_{E^+}: E^+ = E^* \times \R \to M$ and it is clear that
$(\alpha, 0)$ is also a $1$-cocycle of $\tau_{\tilde{E}}:
\tilde{E} = E \times \R \to M$. In addition, if $f(h, \alpha)$ is
the real function on $M$ characterized by
\[
\alpha - h \circ \mu \circ \alpha = f(h, \alpha)1_{E},
\]
it is easy to prove that $f(h, \alpha) = H \circ \alpha$. Thus,
the equation,
\[
d^E(f(h, \alpha)) = 0
\]
is just the Hamilton-Jacobi equation considered in \cite{LMM}.

\subsection{The classical Hamilton-Jacobi equation for time-dependent
Mechanics}\label{ex5.1} Let $\tau:M\to\R$ be a fibration and
$\tau_{1,0}:J^1\tau\to M$ be the associated Lie affgebroid
modelled on the vector bundle $\pi=(\pi_M)_{|V\tau}:V\tau\to M$.
As we know, the bidual vector bundle $\widetilde{J^1\tau}$ to the
affine bundle $\tau_{1,0}:J^1\tau\to M$ may be identified with the
tangent bundle $TM$ to $M$ and, under this identification, the Lie
algebroid structure on $\pi_M:TM\to M$ is the standard Lie
algebroid structure and the $1$-cocycle $1_{J^1\tau}$ on
$\pi_M:TM\to M$ is just the $1$-form $\eta=\tau ^\ast (dt)$, $t$
being the coordinate on $\R$ (see Section \ref{sec1.2}). If
$(t,q^i)$ are local fibred  coordinates on $M$ then $\{
\frac{\partial}{\partial q^i} \}$ (respectively, $\{
\frac{\partial}{\partial t}, \frac{\partial}{\partial q^i} \}$) is
a local basis of sections of $\pi:V\tau\to M$ (respectively,
$\pi_M:TM\to M$). Denote by $(t,q^i,\dot{q}^i)$ (respectively,
$(t,q^i,\dot{t},\dot{q}^i)$) the corresponding local coordinates
on $V\tau$ (respectively, $TM$). Then,  the (local) structure
functions of $TM$ with respect to this local trivialization are
given by
\begin{equation}\label{str-const-trivial}
\begin{array}{l}
C^k_{ij}=0 \mbox{ and } \rho ^i_j=\delta _{ij}, \mbox{ for }i,j, k
\in \{ 0,1,\ldots ,n\}.
\end{array}
\end{equation}

Now, let $\pi^*:V^*\tau\to M$ be the dual vector bundle to
$\pi:V\tau\to M$ and $(J^1\tau)^+\cong T^*M$ be the cotangent
bundle to $M$. Denote by $(t,q^i,p_i)$ (resp., $(t,q^i,p_t,p_i)$)
the dual coordinates on $V^*\tau$ (resp., $T^*M$) to
$(t,q^i,\dot{q}^i)$ (resp., $(t,q^i,\dot{t},\dot{q}^i)$). Then,
since the anchor map of $\pi_M:TM\to M$ is the identity of $TM$,
it follows that the Lie algebroids $\pi_M^{\pi^*}:{\mathcal
T}^{TM}(V^*\tau)\to V^*\tau$ and $\pi_{V^*\tau}:T(V^*\tau)\to
V^*\tau$ are isomorphic.

Next, let $h$ be a Hamiltonian section, that is,  $h:V^*\tau \to
(J^1\tau)^+\cong T^\ast M$ is a section of the canonical
projection $\mu:(J^1\tau)^+\cong T^*M \to V^*\tau$. $h$ is locally
given by
\[
h(t,q^i,p_i)=(t,q^i,-H(t,q^j,p_j),p_i).
\]
Moreover, the cosymplectic structure $(\Omega_h,\eta)$ on the Lie
algebroid $\pi_M^{\pi^*}: {\mathcal T}^{TM}(V^*\tau)\cong
T(V^*\tau) \to V^*\tau$ is, in this case, the standard
cosymplectic structure $(\Omega_h,\eta)$ on the manifold $V^*\tau$
locally given by (see (\ref{Omegahlocal}) and
(\ref{str-const-trivial}))
\[
\Omega_h=dq^i\wedge dp_i + \frac{\partial H}{\partial
q^i}dq^i\wedge dt + \frac{\partial H}{\partial p_i}dp_i\wedge
dt,\;\;\; \; \eta=dt.
\]
Thus, the Reeb section of $(\Omega_h,\eta)$ is the vector field
$R_h$ on $V^*\tau$ defined by
\[
R_h=\frac{\partial }{\partial t}+\frac{\partial H}{\partial
p_i}\frac{\partial }{\partial q^i}-\frac{\partial H}{\partial
q^i}\frac{\partial }{\partial p_i}.
\]
It is clear that the integral curves of $R_h$ \[ t\mapsto
(t,q^i(t),p_i(t))\] are just the solutions of the classical
time-dependent Hamilton equations for $h$
\[
\frac{dq^i}{dt}=\frac{\partial H}{\partial p_i},\quad
\frac{dp_i}{dt}=-\frac{\partial H}{\partial q^i}.
\]
Furthermore, using Theorem \ref{thHJ}, we deduce the following
result.

\begin{corollary}
Let $\tau:M\to\R$ be a fibration and $\tau_{1,0}:J^1\tau\to M$ be
the associated Lie affgebroid modelled on the vector bundle
$\pi=(\pi_M)_{|V\tau}:V\tau\to M$. Let $h:V^*\tau\to T^*M$ be a
Hamiltonian section and $R_h$ be the Reeb vector field of the
corresponding cosymplectic structure $(\Omega_h,\eta)$ on
$V^*\tau$. Suppose that $\alpha$ is a closed $1$-form on $M$ and
denote by $R_h^\alpha$ the vector field on $M$ given by
$R_h^\alpha=T\pi^*\circ R_h\circ\mu\circ\alpha$. Then, the
following conditions are equivalent:
\begin{enumerate}
\item For every integral curve $t\to c(t)$ of $R_h^\alpha$, the curve
$t\to \mu(\alpha(c(t)))$ on $V^*\tau$ satisfies the Hamilton
equations for $h$.
\item $\alpha$ satisfies the Hamilton-Jacobi equation
$d^{V\tau}(f(h,\alpha))=0$, that is, the function $f(h,\alpha)$ is
constant on the leaves of the vertical bundle to $\tau$.
\end{enumerate}
\end{corollary}

\begin{remark}{\em $i)$ We recall that the
function $f(h,\alpha)$ on $M$ is characterized by the following
condition
$$\alpha-h\circ\mu\circ\alpha=f(h,\alpha)\eta.$$

$ii)$ If the fibers of $\tau$ are connected then the equation
$d^{V\tau}(f(h,\alpha))=0$ holds if and only if the function
$f(h,\alpha)$ is constant on the fibers of $\tau$.
\begin{flushright}$\diamondsuit$\end{flushright}}\end{remark}

Now, suppose that the fibration $\tau$ is trivial, that is,
$M=\R\times P$ and $\tau$ is the canonical projection on the first
factor. Then, the vector bundle $\pi^*:V^*\tau\to M$ is isomorphic
to the product $\R\times T^*P$ (as a vector bundle over
$M=\R\times P$). Thus, a Hamiltonian section
$h:V^*\tau\cong\R\times T^*P\to T^*M\cong(\R\times\R)\times T^*P$
may be identified with a time-dependent Hamiltonian function
$H:\R\times T^*P\to\R$ in such a way that
$$h(t,\beta)=(t,-H(t,\beta),\beta),\makebox[1cm]{for}(t,\beta)\in\R\times
T^*P.$$ Moreover, if $\alpha$ is an exact $1$-form on $M$, that
is, $\alpha=dW$ with $W$ a real function on $M=\R\times P$, then
one may prove that the function $f(h,dW)$ is given by
$$f(h,dW)(t,q)=\displaystyle\frac{\partial W}{\partial
t}_{|(t,q)}+H(t,dW_t(q)),\makebox[1cm]{for}(t,q)\in M=\R\times
P,$$ where $W_t:P\to\R$ is the real function on $P$ defined by
$$W_t(q)=W(t,q).$$
Thus, if $P$ is connected the equation
$$d^{V\tau}(f(h,dW))=0$$
holds if and only if the function
$$(t,q)\in\R\times P\mapsto\displaystyle\frac{\partial W}{\partial
t}_{|(t,q)}+H(t,dW_t(q))\in\R$$ doesn't depend on $q$.

This last condition may be locally expressed as follows:
$$\displaystyle\frac{\partial W}{\partial
t}+H(t,q^i,\displaystyle\frac{\partial W}{\partial
q^i})=\makebox{constant on $P$}$$ which is the classical
time-dependent Hamilton-Jacobi equation for the function $W$ (see
\cite{AM}).

\subsection{The Hamilton-Jacobi equation on the Atiyah
affgebroid}\label{ex5.2} Let $p:Q\to M$ be a principal $G$-bundle.
Denote by $\Phi:G\times Q\to Q$ the free action of $G$ on $Q$ and
by $T\Phi:G\times TQ\to TQ$ the tangent action of $G$ on $TQ.$
Then, one may consider the quotient vector bundle
$\pi_{Q}|G:TQ/G\to M=Q/G$ and the sections of this vector bundle
may be identified with the vector fields on $Q$ which are
invariant under the action $\Phi$. Using that every $G$-invariant
vector field on $Q$ is $p$-projectable and that the usual Lie
bracket on vector fields is closed with respect to $G$-invariant
vector fields, we can induce a Lie algebroid structure on $TQ/G$.
This Lie algebroid is called {\it the Atiyah algebroid} associated
with the principal $G$-bundle $p:Q\to M$ (see \cite{LMM,Ma}).

Now, suppose that $\nu:M\to \R$ is a fibration of $M$ on $\R$.
Denote by $\tau:Q\to \R$ the composition $\tau=\nu\circ p.$ Then,
$\Phi$ induces an action $J^1\Phi:G\times J^1\tau\to J^1\tau$ of
$G$ on $J^1\tau$ such that
$$J^1\Phi(g,j_t^1\gamma)= j_t^1(\Phi_g\circ \gamma),$$ for all $g\in
G$ and $\gamma:I\subset\R\to Q$ a local section of $\tau$  with
$t\in I.$ Moreover, the projection
\[
\tau_{1,0}|G:J^1\tau/G\to M,\;\;\;\; [j^1_t\gamma]\mapsto
p(\tau_{1,0}(j^1_t\gamma))=p(\gamma(t))
\]
defines an affine bundle on $M$ which is modelled on the quotient
vector bundle
\[
\pi|G:V\tau/G\to M,\;\;\; [u_q]\mapsto p(q), \mbox{ for } u_q\in
V_q\tau,
\]
$\pi:V\tau\to Q$ being the vertical bundle of the fibration
$\tau:Q\to \R.$ Here, the action of $G$ on $V\tau$ is the
restriction to $V\tau$ of the tangent action $T\Phi$ of $G$ on
$TQ.$

In addition, the bidual vector bundle of $J^1\tau/G\to M$ is
$\pi_{Q}|G:TQ/G\to M$.

On the other hand, if $t$ is the usual coordinate on $\R$, the
$1$-form $\eta=\tau^*(dt)$ is $G$-invariant and defines a non-zero
$1$-cocycle $\phi:TQ/G\to \R$ on the Atiyah algebroid $TQ/G$. Note
that $\phi^{-1}\{1\}\cong J^1\tau/G$ and therefore, one may
consider the corresponding Lie affgebroid structure on $J^1\tau/G$
(see \cite{MMeS}). $J^1\tau/G$ endowed with this structure is
called {\it the Atiyah affgebroid }  associated with the principal
$G$-bundle $p:Q\to M$ and the fibration $\nu:M\to \R$ (see
\cite{IMPS,MMeS}).

The Lie group $G$ acts on the vector bundles $T^*Q$ and $V^*\tau$
in such a way that the dual vector bundles to $TQ/G$ and $V\tau/G$
may be identified with the quotient vector bundles $T^*Q/G$ and
$V^*\tau/G$, respectively (see \cite{LMM}). Moreover, the
canonical projection between the vector bundles $(J^1\tau/G)^+
\cong T^*Q/G$ and $V^*\tau/G$ is the map $\mu|G$ given by
\[
(\mu|G)[\alpha_{q}] = [\mu(\alpha_{q})], \makebox[.3cm]{} \mbox{
for } \alpha_{q} \in T_{q}^*Q \mbox{ and } q \in Q,
\]
where $\mu: T^*Q \to V^*\tau$ is the projection between $T^*Q$ and
$V^*\tau$. Thus, if $\tilde{h}: V^*\tau/G \to (J^1\tau/G)^+\cong
T^*Q/G$ is a Hamiltonian section of $\mu|G$ then $\tilde{h}$
induces a $G$-equivariant Hamil\-to\-nian section $h: V^*\tau \to
T^*Q$ such that $\tilde{h} = h|G$, that is,
\[
\tilde{h}[\alpha_{q}] = [h(\alpha_{q})], \makebox[.3cm]{} \mbox{
for } \alpha_{q} \in V_{q}^*\tau \mbox{ and } q \in Q.
\]
Conversely, if $h: V^*\tau \to T^*Q$ is a $G$-equivariant
Hamiltonian section of $\mu: T^*Q \to V^*\tau$ then $h$ induces a
Hamiltonian section $\tilde{h}: V^*\tau/G \to T^*Q/G$ such that
$\tilde{h} = h|G$, that is,
\[
\tilde{h}[\alpha_{q}] = [h(\alpha_{q})], \makebox[.3cm]{} \mbox{
for } \alpha_{q} \in V_{q}^*\tau \mbox{ and } q\in Q.
\]
Next, we will discuss the relation between the solutions of the
Hamilton-Jacobi equation for the Hamiltonians $h$ and $\tilde{h}$.
In fact, we will prove the following result.
\begin{proposition}
There exists a one-to-one correspondence between the solutions of
the Ha\-mil\-ton-Jacobi equation for $h|G$ and the $G$-invariant
solutions of the Hamilton-Jacobi equation for $h$.
\end{proposition}
\begin{proof}
If $p_{TQ}: TQ \to TQ/G$ is the canonical projection then $p_{TQ}$
is a fiberwise bijective Lie algebroid morphism over $p: Q \to
M=Q/G$. Thus, there exists a one-to-one correspondence between the
$1$-cocycles of the Atiyah algebroid $\tau_{Q}|G: TQ/G \to M= Q/G$
and the $G$-invariant closed $1$-forms on $Q$. Indeed, if $\alpha:
Q \to T^*Q$ is a $G$-invariant closed $1$-form on $Q$ then
$\alpha|G: M= Q/G \to T^*Q/G$ defined by
\[
(\alpha|G)([q]) = [\alpha(q)], \makebox[.3cm]{} \mbox{ for } q \in
Q,
\]
is a $1$-cocycle of the Atiyah algebroid (note that $\alpha|G =
(p_{TQ}, p)^*(\alpha)$). Moreover, if
\[
\alpha - h \circ \mu \circ \alpha = f(h, \alpha)\eta
\]
then
\[
\alpha|G - (h|G) \circ (\mu|G) \circ (\alpha|G) = (f(h,
\alpha)|G)\phi,
\]
where $f(h, \alpha)|G: M = Q/G \to \R$ is the real function on $M$
which is characterized by the condition
\begin{equation}\label{funciones}
(f(h, \alpha)|G) \circ p = f(h, \alpha).
\end{equation}
Therefore, $f(h|G, \alpha|G) = f(h, \alpha)|G$.

On the other hand, if $p_{V\tau}: V\tau \to V\tau/G$ is the
canonical projection then $p_{V\tau}$ is a Lie algebroid morphism
over $p: Q \to M=Q/G$ and, from (\ref{funciones}), it follows that
\[
(p_{V\tau}, p)^*(d^{V\tau/G}(f(h, \alpha)|G)) = d^{V\tau}(f(h,
\alpha)).
\]
Finally, using that $p_{V\tau}$ is a fiberwise bijective morphism,
we conclude that
\[
d^{V\tau/G}(f(h,\alpha)|G)=0 \Leftrightarrow
d^{V\tau}(f(h,\alpha))=0,
\]
which proves the result.
\end{proof}

\end{document}